\documentclass{amsart}%
\usepackage{graphicx}
\usepackage{amscd}
\usepackage{amsmath}
\usepackage{amsfonts}
\usepackage{amssymb}%
\newtheorem{theorem}{Theorem}
\theoremstyle{plain}

\newtheorem{example}{Example}

\newtheorem{lemma}{Lemma}

\numberwithin{equation}{section}

\begin{document}
\title[Banach spaces whose duals are isomorphic to $l_{1}(\Gamma)$]{On Banach spaces whose duals are isomorphic to $l_{1}(\Gamma)$}
\author{Daniel M. Pellegrino}
\address[Daniel M. Pellegrino]{Depto de Matem\'{a}tica e Estat\'{\i}stica- Caixa Postal
10044- UFCG- Campina Grande-PB-Brazil }
\email{dmp@dme.ufcg.edu.br}
\thanks{This work is partially supported by Instituto do Mil\^{e}nio, IMPA}

\begin{abstract}
In this paper we present new characterizations of Banach spaces whose duals
are isomorphic to $l_{1}(\Gamma),$ extending results of Stegall, Lewis-Stegall
and Cilia-D'Anna-Guti\'{e}rrez.

\end{abstract}
\maketitle

\section{Introduction, notation and background}

Banach spaces whose duals are isomorphic to $l_{1}(\Gamma)$ were investigated
in the works of Stegall \cite{Stegall}, Lewis-Stegall \cite{Lewis} and, in a
recent paper, Cilia-D'Anna-Guti\'{e}rrez \cite{Cilia} studied polynomial
characterizations of such spaces. Our aim is to show that polynomial
characterizations of Banach spaces whose duals are isomorphic to $l_{1}%
(\Gamma)$ are extremely much more common than it is known, and many of these
statements are consequences of simple results concerning polynomial ideals.
Our techniques also give alternative non-tensorial proofs for the results of
\cite{Cilia}.

Throughout this paper $E,E_{1},...,E_{n},F,G,$ will stand for (real or
complex) Banach spaces, $B_{E}$ will denote the closed unit ball on $E$ and
$\mathbb{N}$ will denote the set of the natural numbers.

The space of all continuous $n$-homogeneous polynomials from $E$ into $F$
endowed with the $\sup$ norm is represented by $\mathcal{P}(^{n}E;F)$ and the
space of all continuous $n$-linear mappings from $E_{1}\times...\times E_{n}$
into $F$ (with the $\sup$ norm) is denoted by $\mathcal{L}(E_{1}%
,...,E_{n};F).$ When $E_{1}=...=E_{n}=E$, we write $\mathcal{L}(^{n}E;F).$ If
$P\in\mathcal{P}(^{n}E;F),$ we use the symbol $\overset{\vee}{P}$ for the
(unique) symmetric $n$-linear mapping associated to $P$. On the other hand, if
$T\in\mathcal{L}(^{n}E;F)$ we write $\overset{\wedge}{T}(x)=T(x,...,x).$ For
$i=1,...,n,$ $\Psi_{i}^{(n)}:\mathcal{L}(E_{1},...,E_{n};F)\rightarrow
\mathcal{L}(E_{i};\mathcal{L}(E_{1},\overset{[i]}{...},E_{n};F))$ will
represent the canonical isometric isomorphism given by
\[
\Psi_{i}^{(n)}(T)(x_{i})(x_{1}\overset{[i]}{...}x_{n})=T(x_{1},...,x_{n}),
\]
where the notation $\overset{[i]}{...}$ means that the $i$-th coordinate is
not involved.

An ideal of (homogeneous) polynomials $\mathfrak{P}$ is a subclass of the
class of all continuous homogeneous polynomials between Banach spaces such
that for all index $n$ and all $E$ and $F$, \ the components $\mathfrak{P}%
(^{n}E;F)=\mathcal{P}(^{n}E;F)\cap\mathfrak{P}$ satisfy:

(i) $\mathfrak{P}(^{n}E;F)$ is a linear subspace of $\mathcal{P}(^{n}E;F)$
which contains the polynomials of finite type.

(ii) If $P\in\mathfrak{P}(^{n}E;F),$ $T_{1}\in\mathcal{L}(G;E)$ and $T_{2}%
\in\mathcal{L}(F;H),$ then $T_{2}PT_{1}\in\mathfrak{P}(^{n}G;H).$

In this note we will be concerned with two special methods of creating ideals
of polynomials: factorization and linearization.

\begin{itemize}
\item (Factorization method) If $\mathfrak{I}$ is an operator ideal, a
continuous $n$-homogeneous polynomial $P\in\mathcal{P}(^{n}E;F)$ is of type
$\mathcal{P}_{\mathcal{L}(\mathfrak{I})}$ if there exists a Banach space $G$,
a linear operator $T\in\mathfrak{I}(E;G)$ and $Q\in\mathcal{P}(^{n}G;F)$ such
that $P=QT.$

\item (Linearization method) If $\mathfrak{I}$ is an operator ideal,
$T\in\mathcal{L}(E_{1},...,E_{n};F)$ is of type $[\mathfrak{I]}$ if $\Psi
_{i}^{(n)}(T)\in\mathfrak{I}(E_{i};\mathcal{L}(E_{1},\overset{[i]}{...}%
,E_{n}))$ for every $i=1,...,n.$ We say that $P\in\mathcal{P}(^{n}E;F)$ is of
type $\mathcal{P}_{[\mathfrak{I]}}$ if $\overset{\vee}{P}$ is of type
$[\mathfrak{I]}$.
\end{itemize}

An $n$-homogeneous polynomial is said to be $p$-dominated if there exist
$C\geq0$ and a regular probability measure $\mu$ \ on the Borel $\sigma
$-algebra on B$_{E^{\prime}}$(with the weak star topology) such that
\begin{equation}
\left\Vert P\left(  x\right)  \right\Vert \leq C\left[  \int_{B_{E^{\prime}}%
}\left\vert \varphi\left(  x\right)  \right\vert ^{p}d\mu\left(
\varphi\right)  \right]  ^{\frac{n}{p}}.
\end{equation}

We write $\mathcal{P}_{d,p}(^{n}E;F)$ to denote the space of $p$-dominated
$n$-homogeneous polynomials from $E$ into $F.$ For $n=1$ we obtain the
$p$-absolutely summing operator. We represent the space of all absolutely
$p$-summing operators from $E$ into $F$ by $\mathcal{L}_{as,p}(E;F).$ It is
well known that $\mathcal{P}_{d,p}(^{n}E;F)=\mathcal{P}_{\mathcal{L}%
(as,p)}(^{n}E;F)$. For references on $p$-dominated polynomials we mention
(\cite{irish},\cite{Matos2},\cite{studia}, among others). For details
concerning polynomials on Banach spaces we mention \cite{Dineen}.

\section{Results}

We shall start with some useful Lemmas:

\begin{lemma}
\label{aaa} If $\mathfrak{I}_{1}$ and $\mathfrak{I}_{2}$ are ideals of
polynomials, and
\begin{equation}
\mathcal{L}_{\mathfrak{I}_{1}}(E;F)\subset\mathcal{L}_{\mathfrak{I}_{2}%
}(E;F)\text{ for every }F \label{nova}%
\end{equation}
then
\[
\mathcal{P}_{\mathcal{L}[\mathfrak{I}_{1}]}(^{m}E;F)\subset\mathcal{P}%
_{\mathcal{L}[\mathfrak{I}_{2}]}(^{m}E;F)\text{ and }\mathcal{P}%
_{[\mathfrak{I}_{1}]}(^{m}E;F)\subset\mathcal{P}_{[\mathfrak{I}_{2}]}%
(^{m}E;F)\text{ for every }F.
\]

\end{lemma}

Proof. If $P\in\mathcal{P}_{\mathcal{L}[\mathfrak{I}_{1}]}(^{m}E;F),$ then
$P=Qu,$ with $Q\in\mathcal{P}(^{m}G;F)$ and $u\in\mathcal{L}_{\mathfrak{I}%
_{1}}(E;G).$ From (\ref{nova}), we have $u\in\mathcal{L}_{\mathfrak{I}_{2}%
}(E;G)$ and thus $P\in\mathcal{P}_{\mathcal{L}[\mathfrak{I}_{2}]}(^{m}E;F).$
The other case is similar.$\Box$

\begin{lemma}
\label{teoee}If $\mathcal{I}_{1}$ and $\mathcal{I}_{2}$ are ideals of
polynomials and $\mathcal{P}_{\mathcal{I}_{0}}(^{n}E;F)\cap\mathcal{P}%
_{\mathcal{I}_{1}}(^{n}E;F)\subset\mathcal{P}_{\mathcal{I}_{2}}(^{n}E;F)$ for
some natural $n,$ suppose that the following hold true:
\end{lemma}

(i) $P\in\mathcal{P}_{\mathcal{I}_{2}}(^{n}E;F)\Rightarrow\overset{\vee}%
{P}(.,a,...,a)\in$ $\mathcal{L}_{\mathcal{I}_{2}}(E;F)$ for every $a\in E$, fixed.

(ii) For $j=0,1$ and $m<n,$ if $P\in\mathcal{P}_{\mathcal{I}_{j}}(^{m}E;F)$
and $\varphi\in\mathcal{L}(E;\mathbb{K}),$ then $P.\varphi\in\mathcal{P}%
_{\mathcal{I}_{j}}(^{m+1}E;F).$

Then $\mathcal{L}_{\mathcal{I}_{0}}(E;F)\cap\mathcal{L}_{\mathcal{I}_{1}%
}(E;F)\subset\mathcal{L}_{\mathcal{I}_{2}}(E;F).$

Proof. If $T\in\mathcal{L}_{\mathcal{I}_{0}}(E;F)\cap\mathcal{L}%
_{\mathcal{I}_{1}}(E;F)$, then define $\varphi\in$ $\mathcal{L}(E;\mathbb{K}%
),$ $\varphi\neq0$ and $a\in E$ such that $\varphi(a)=1$ and consider the
following $n$-homogeneous polynomial:
\[
P(x)=T(x)\varphi(x)^{n-1}.
\]
By applying (ii), $P\in\mathcal{P}_{\mathcal{I}_{0}}(^{n}E;F)\cap
\mathcal{P}_{\mathcal{I}_{1}}(^{n}E;F)\subset\mathcal{P}_{\mathcal{I}_{2}%
}(^{n}E;F).$ Finally (i) yields that $\overset{\vee}{P}(.,a,...,a)\in
\mathcal{L}_{\mathcal{I}_{2}}(E;F)$ and thus
\[
\frac{1}{n}T+\frac{n-1}{n}T(a)\varphi\in\mathcal{L}_{\mathcal{I}_{2}}(E;F).
\]
Since $\varphi\in\mathcal{L}_{\mathcal{I}_{2}}(E;\mathbb{K}),$ we conclude
that $T\in\mathcal{L}_{\mathcal{I}_{2}}(E;F).\Box$

\begin{lemma}
\label{novo}If $\mathcal{P}_{\mathcal{L}[\mathcal{I}_{0}]}(^{n}E;F)\cap
\mathcal{P}_{\mathcal{L}[\mathcal{I}_{1}]}(^{n}E;F)\subset\mathcal{P}%
_{\mathcal{I}_{2}}(^{n}E;F)$ and $\mathcal{P}_{\mathcal{I}_{2}}(^{n}E;F)$
satisfies the hypothesis (i) of Lemma \ref{teoee}, then%
\[
\mathcal{L}_{\mathcal{I}_{0}}(E;F)\cap\mathcal{L}_{\mathcal{I}_{1}%
}(E;F)\subset\mathcal{L}_{\mathcal{I}_{2}}(E;F).
\]

\end{lemma}

Proof. If $T\in\mathcal{L}_{\mathcal{I}_{0}}(E;F)\cap\mathcal{L}%
_{\mathcal{I}_{1}}(E;F),$ choosing a continuous (non null) linear functional
$\varphi$ on $F,$ define an $n$-homogeneous polynomial $P:E\rightarrow F$ by
$P(x)=T(x)\varphi^{n-1}(T(x)).$ Then $P=Q\circ T$, where $Q:F\rightarrow F$ is
given by $Q(y)=y\varphi^{n-1}(y).$ Thus $P\in\mathcal{P}_{\mathcal{L}%
[\mathcal{I}_{0}]}(^{n}E;F)\cap\mathcal{P}_{\mathcal{L}[\mathcal{I}_{1}]}%
(^{n}E;F)\subset\mathcal{P}_{\mathcal{I}_{2}}(^{n}E;F)$ and since
$\mathcal{P}_{\mathcal{I}_{2}}(^{n}E;F)$ satisfies the hypothesis (i) of Lemma
\ref{teoee}, we have that $\overset{\vee}{P}(.,a,...,a)\in$ $\mathcal{L}%
_{\mathcal{I}_{2}}(E;F)$ (for $a\in E$ so that $\varphi(a)\neq0$) and hence
$T\in\mathcal{L}_{\mathcal{I}_{2}}(E;F).\Box$

Let us recall the concepts of compact and nuclear polynomials. A polynomial
$P:E\rightarrow F$ is said to be compact if $P(B_{E})$ is relatively compact
in $F$. The space of all compact $m$-homogeneous polynomials from $E$ into $F$
will be denoted by $\mathcal{P}_{K}(^{m}E;F).$ For the compact operators from
$E$ into $F$ we use the symbol $\mathcal{L}_{K}(E;F).$ We say that
$P\in\mathcal{P}(^{m}E;F)$ is nuclear if it is possible to find $(\varphi
_{i})\subset E^{\prime}$ and $(y_{i})\subset F$ so that
\[
Px=\sum\limits_{i=1}^{\infty}\left[  \varphi_{i}(x)\right]  ^{m}y_{i}\text{
and }\sum\limits_{i=1}^{\infty}\left\|  \varphi_{i}\right\|  ^{m}\left\|
y_{i}\right\|  <\infty.
\]
The space of all nuclear $m$-homogeneous polynomials from $E$ into $F$ is
denoted by $\mathcal{P}_{N}(^{m}E;F).$ For the linear case we write
$\mathcal{L}_{N}(E;F).$ The relation between nuclear, compact operators
(polynomials), and Banach spaces whose duals are isomorphic to $l_{1}(\Gamma)$
is given by the following results:

\begin{theorem}
\label{T2}(Lewis-Stegall \cite{Lewis}/Stegall \cite{Stegall}) Given a Banach
space $E$, the following assertions are equivalent:

(i) $E^{\prime}$ is isomorphic to $l_{1}(\Gamma)$ for some $\Gamma.$

(ii) For every Banach space $F$, we have $\mathcal{L}_{as,1}(E;F)\subset
\mathcal{L}_{N}(E;F).$

(iii) For every Banach space $F$, $\mathcal{L}_{as,1}(E;F)\cap\mathcal{L}%
_{K}(E;F)\subset\mathcal{L}_{N}(E;F).$
\end{theorem}

\begin{theorem}
\label{T}(Cilia-D'Anna-Guti\'{e}rrez \cite{Cilia}) Given a Banach space $E$,
the following assertions are equivalent:

(i) $E^{\prime}$ is isomorphic to $l_{1}(\Gamma)$ for some $\Gamma.$

(ii) For all natural $m$ and every Banach space $F$, we have $\mathcal{P}%
_{d,1}(^{m}E;F)\subset\mathcal{P}_{N}(^{m}E;F).$

(iii) There is a natural $m$ such that for every $F$ we have $\mathcal{P}%
_{d,1}(^{m}E;F)\subset\mathcal{P}_{N}(^{m}E;F).$

(iv) There is a natural $m$ such that for every $F$ we have $\mathcal{P}%
_{d,1}(^{m}E;F)\cap\mathcal{P}_{K}(^{m}E;F)\subset\mathcal{P}_{N}(^{m}E;F).$
\end{theorem}

Our results will show that it is possible to show a considerably longer list
of characterizations of Banach spaces whose duals are isomorphic to
$l_{1}(\Gamma)$ and present different proofs for each assertion of Theorem
\ref{T}.

\begin{theorem}
Given a Banach space $E$, the following assertions are equivalent:

(i) $E^{\prime}$ is isomorphic to $l_{1}(\Gamma)$ for some $\Gamma.$

(ii) For all $m\in\mathbb{N}$ and every $F$, we have $\mathcal{P}_{d,1}%
(^{m}E;F)\subset\mathcal{P}_{\mathcal{L}[N]}(^{m}E;F).$

(iii) There is $m\in\mathbb{N}$ such that for every $F$ we have $\mathcal{P}%
_{d,1}(^{m}E;F)\subset\mathcal{P}_{\mathcal{L}[N]}(^{m}E;F).$

(iv) For all $m\in\mathbb{N}$ and every $F$, we have $\mathcal{P}_{d,1}%
(^{m}E;F)\subset\mathcal{P}_{N}(^{m}E;F).$

(v) There is $m\in\mathbb{N}$ such that for every $F$ we have $\mathcal{P}%
_{d,1}(^{m}E;F)\subset\mathcal{P}_{N}(^{m}E;F).$

(vi) For all $m\in\mathbb{N}$ and every $F$, we have $\mathcal{P}_{d,1}%
(^{m}E;F)\subset\mathcal{P}_{[N]}(^{m}E;F).$

(vii) There is $m\in\mathbb{N}$ such that for every $F$ we have $\mathcal{P}%
_{d,1}(^{m}E;F)\subset\mathcal{P}_{[N]}(^{m}E;F).$

(viii) For all $m\in\mathbb{N}$ and every $F$, we have $\mathcal{P}%
_{[as,1]}(^{m}E;F)\subset\mathcal{P}_{[N]}(^{m}E;F).$

(ix) There is $m\in\mathbb{N}$ such that for every $F$ we have $\mathcal{P}%
_{[as,1]}(^{m}E;F)\subset\mathcal{P}_{[N]}(^{m}E;F).$

(x) For all $m\in\mathbb{N}$ and every $F$ we have $\mathcal{P}_{d,1}%
(^{m}E;F)\cap\mathcal{P}_{K}(^{m}E;F)\subset\mathcal{P}_{N}(^{m}E;F).$

(xi) There is $m\in\mathbb{N}$ such that for every $F$ we have $\mathcal{P}%
_{d,1}(^{m}E;F)\cap\mathcal{P}_{K}(^{m}E;F)\subset\mathcal{P}_{N}(^{m}E;F).$

(xii) For all $m\in\mathbb{N}$ and every $F$, we have $\mathcal{P}%
_{[as,1]}(^{m}E;F)\cap\mathcal{P}_{[K]}(^{m}E;F)\subset\mathcal{P}_{[N]}%
(^{m}E;F).$

(xiii) There is $m\in\mathbb{N}$ such that for every $F$ we have
$\mathcal{P}_{[as,1]}(^{m}E;F)\cap\mathcal{P}_{[K]}(^{m}E;F)\subset
\mathcal{P}_{[N]}(^{m}E;F).$

(xiv) For all $m\in\mathbb{N}$ and every $F$, we have $\mathcal{P}_{d,1}%
(^{m}E;F)\cap\mathcal{P}_{[K]}(^{m}E;F)\subset\mathcal{P}_{\mathcal{L}%
[N]}(^{m}E;F).$

(xv) There is $m\in\mathbb{N}$ and every $F$, we have $\mathcal{P}_{d,1}%
(^{m}E;F)\cap\mathcal{P}_{[K]}(^{m}E;F)\subset\mathcal{P}_{\mathcal{L}%
[N]}(^{m}E;F).$

(xvi) For all $m\in\mathbb{N}$ such that for every $F$ we have $\mathcal{P}%
_{d,1}(^{m}E;F)\cap\mathcal{P}_{[K]}(^{m}E;F)\subset\mathcal{P}_{[N]}%
(^{m}E;F).$

(xvii) There is $m\in\mathbb{N}$ such that for every $F$ we have
$\mathcal{P}_{d,1}(^{m}E;F)\cap\mathcal{P}_{[K]}(^{m}E;F)\subset
\mathcal{P}_{[N]}(^{m}E;F).$
\end{theorem}

Proof.

(i)$\Rightarrow$(ii) is consequence of the Theorem of Lewis-Stegall and Lemma
\ref{aaa}. (ii)$\Rightarrow$(iii) is obvious.

A direct computation gives $\mathcal{P}_{\mathcal{L}[N]}(^{m}E;F)\subset
\mathcal{P}_{N}(^{m}E;F)$ and hence it is easy to see that (iii)$\Rightarrow
$(iv)$\Rightarrow$(v). It is not hard to check that the ideals of nuclear
polynomials and dominated polynomials satisfy (i) and (ii) of Lemma
\ref{teoee}, respectively (these facts will be used several times in the
present proof). Hence (v) implies $\mathcal{L}_{as,1}(E;F)\subset
\mathcal{L}_{N}(E;F)$ and consequently we obtain (i).

(ii)$\Rightarrow$(vi) holds because $\mathcal{P}_{\mathcal{L}[N]}%
(^{m}E;F)\subset\mathcal{P}_{[N]}(^{m}E;F)$ (it is true for arbitrary ideals
of polynomials)$.$ (vi)$\Rightarrow$(vii) is obvious.

In order to prove (vii)$\Rightarrow$(i) it suffices to show that (vii) implies
$\mathcal{L}_{as,1}(E;F)\subset\mathcal{L}_{N}(E;F).$ So, in order to apply
Lemma \ref{teoee} we must show that whenever $P\in\mathcal{P}_{[N]}(^{m}E;F)$
we have $\overset{\vee}{P}(.,a,...,a)\in\mathcal{L}_{N}(E;F).$ In fact, if
$P\in\mathcal{P}_{[N]}(^{m}E;F)$, we can find $(\varphi_{i})\subset E^{\prime
}$ and $(y_{i})\subset\mathcal{L}(^{m-1}E;F)$ so that
\[
\Psi_{1}^{(m)}(\overset{\vee}{P})(x)=\sum\limits_{i=1}^{\infty}\left[
\varphi_{i}(x)\right]  y_{i}\text{ and }\sum\limits_{i=1}^{\infty}\left\Vert
\varphi_{i}\right\Vert \left\Vert y_{i}\right\Vert <\infty.
\]
Thus
\[
\overset{\vee}{P}(x,a,...,a)=\Psi_{1}^{(m)}(\overset{\vee}{P}%
)(x)(a,...,a)=\sum\limits_{i=1}^{\infty}\left[  \varphi_{i}(x)\right]
y_{i}(a,...,a)
\]
and
\[
\sum\limits_{i=1}^{\infty}\left\Vert \varphi_{i}\right\Vert \left\Vert
y_{i}(a,...,a)\right\Vert \leq\sum\limits_{i=1}^{\infty}\left\Vert \varphi
_{i}\right\Vert \left\Vert y_{i}\right\Vert \left\Vert a\right\Vert
^{m-1}<\infty.
\]
Hence $\mathcal{P}_{[N]}(^{m}E;F)$ satisfy (i) of \ Lemma \ref{teoee}.

(i)$\Rightarrow$(viii) is due to the result of Lewis-Stegall and Lemma
\ref{aaa}. (viii)$\Rightarrow$(ix) is obvious. For the proof of
(ix)$\Rightarrow$(i) one may realize that a standard use of Ky Fan's Lemma
yields that a continuous (symmetric) multilinear mapping $T:E\times...\times
E\rightarrow F$ is of type $[as,p]$ if, and only if,\ there exist $C\geq0$ and
a regular probability measure $\mu\in P\left(  B_{E^{\prime}}\right)  ,$ such
that
\begin{equation}
\left\Vert T\left(  x_{1},...,x_{n}\right)  \right\Vert \leq C\left\Vert
x_{1}\right\Vert ...\left\Vert x_{n-1}\right\Vert \left[  \int_{B_{E^{\prime}%
}}\left\vert \varphi\left(  x_{n}\right)  \right\vert ^{p}d\mu\left(
\varphi\right)  \right]  ^{\frac{1}{p}}%
\end{equation}
and thus $\mathcal{P}_{[as,1]}(^{m}E;F)$ satisfy (ii) of Lemma \ref{teoee} and
we conclude that (ix) implies $\mathcal{L}_{as,1}(E;F)\subset\mathcal{L}%
_{N}(E;F).$ Thus, the result of Lewis-Stegall completes the proof.

(iv)$\Rightarrow$(x)$\Rightarrow$(xi) is trivial. Since the ideal of compact
polynomials also satisfies (ii) of Lemma \ref{teoee}, (xi) implies
$\mathcal{L}_{as,1}(E;F)\cap\mathcal{L}_{K}(E;F)\subset\mathcal{L}_{N}(E;F)$
for every $F$ and thus we obtain (i).

In order to prove (i)$\Rightarrow$(xii) we observe that (i) implies
$\mathcal{L}_{as,1}(E;F)\cap\mathcal{L}_{K}(E;F)\subset\mathcal{L}_{N}(E;F)$
for every $F$ and thus Lemma \ref{aaa} asserts that (xii) holds.

(xii)$\Rightarrow$(xiii) is obvious.

In order to prove (xiii)$\Rightarrow$(i) we must show that $\mathcal{P}%
_{[K]}(^{m}E;F)$ satisfy (ii) of Lemma \ref{teoee}. If $P\in\mathcal{P}%
_{[K]}(^{m}E;F)$ and $\varphi\in\mathcal{L}(E;F)$, we shall firstly prove that
$R$ defined by $R(x_{1},...,x_{m+1})=\frac{1}{m+1}\varphi(x_{1})\overset{\vee
}{P}(x_{2},...,x_{m+1})...+\frac{1}{m+1}\varphi(x_{m+1})\overset{\vee}%
{P}(x_{1},...,x_{m})$ is of type $[K].$ In fact, since
\[
\Psi_{1}^{(m+1)}(R)(x)=\frac{1}{m+1}\varphi(x).\overset{\vee}{P}+\frac{m}%
{m+1}\varphi.\Psi_{1}^{(m)}(\overset{\vee}{P})(x),
\]
and since $\varphi$ and $\Psi_{1}^{(m)}(\overset{\vee}{P})$ are compact
mappings, we conclude that $\Psi_{1}^{(m+1)}(R)$ is compact. Thus $R$ is of
type $[K]$ and hence $\varphi.P\in\mathcal{P}_{[K]}(^{m+1}E;F).\ $So,
$\mathcal{P}_{[K]}(^{m}E;F)$ satisfy (ii) of Lemma \ref{teoee} and we conclude
that $\mathcal{L}_{as,1}(E;F)\cap\mathcal{L}_{K}(E;F)\subset\mathcal{L}%
_{N}(E;F)$ and obtain (i).

Since the ideal of compact operators is closed, injective and surjective, we
have that $\mathcal{P}_{[K]}=\mathcal{P}_{\mathcal{L}[K]}$ and this fact will
be used in each one of the next arguments$.$ For the proof of (i)$\Rightarrow
$(xiv) note that (i) implies $\mathcal{L}_{as,1}(E;F)\cap\mathcal{L}%
_{K}(E;F)\subset\mathcal{L}_{N}(E;F)$ and Lemma \ref{aaa} furnishes the proof.
The proof that (xiv) implies (xv) is immediate. Since $\mathcal{P}%
_{\mathcal{L}[N]}(E;F)\subset\mathcal{P}_{[N]}(E;F)$ we obtain
(xv)$\Rightarrow$(xvi). Finally, (xvi)$\Rightarrow$(xvii) is trivial and
(xvii)$\Rightarrow$(i) is obtained by invoking Lemma \ref{novo} .\bigskip
$\Box$

It is worth remarking, for example, that $\mathcal{P}_{d,1}(^{n}E;F)$ and
$\mathcal{P}_{[as,1]}(^{n}E;F)$ are different spaces, in general, showing that
our results are different from the previous characterizations given in
Theorems \ref{T2} and \ref{T}. The following example was suggested by Prof. M.
C. Matos.

\begin{example}
Define $P:l_{2}\rightarrow\mathbb{\ K}$ by $P(x)=%
{\displaystyle\sum\limits_{j=1}^{\infty}}
\frac{1}{j^{\alpha}}x_{j}^{2}$ with $\alpha=\frac{1}{2}+\varepsilon$ and
$0<\varepsilon<\frac{1}{2}.$ Then $P\in\mathcal{P}_{[as,1]}(^{2}%
l_{2};\mathbb{K})$ and $P\notin\mathcal{P}_{d,1}(^{2}l_{2};\mathbb{K})$ .
\end{example}

In fact, $\overset{\vee}{P}:l_{2}\times l_{2}\rightarrow\mathbb{K}\ $\ is
given by $\overset{\vee}{P}(x,y)=\sum\limits_{j=1}^{\infty}\frac{1}{j^{\alpha
}}x_{j}y_{j}$ and $(\frac{1}{j^{\alpha}})_{j=1}^{\infty}\in l_{2}.$

It suffices to show that $\overset{\vee}{P}$ fails to be $1$-dominated, and
$\Psi_{1}^{(2)}(\overset{\vee}{P})\in\mathcal{L}_{as,1}(l_{2};l_{2}).$ Since
\[
\left(  \sum\limits_{j=1}^{m}\left\Vert \overset{\vee}{P}(e_{j},e_{j}%
)\right\Vert ^{\frac{1}{2}}\right)  ^{2}=\left[  \sum\limits_{j=1}^{m}\left(
\frac{1}{j^{\alpha}}\right)  ^{\frac{1}{2}}\right]  ^{2}\geq\left[
\sum\limits_{j=1}^{m}\left(  \frac{1}{m^{\frac{\alpha}{2}}}\right)  \right]
^{2}=m^{2-\alpha},
\]
if we had
\[
\left(  \sum\limits_{j=1}^{m}\left\Vert \overset{\vee}{P}(e_{j},e_{j}%
)\right\Vert ^{\frac{1}{2}}\right)  ^{2}\leq C\left\Vert (e_{j})_{j=1}%
^{m}\right\Vert _{w,1}^{2},
\]
we would obtain $m^{2-\alpha}\leq C(m^{\frac{1}{2}})^{2}=Cm $ and it is a
contradiction since $\alpha<1.$

In order to prove that $\Psi_{1}^{(2)}(\overset{\vee}{P})\in\mathcal{L}%
_{as,1}(l_{2};l_{2}),$ observe that
\[
\Psi_{1}^{(2)}(\overset{\vee}{P})((x_{j})_{j=1}^{\infty})=\left(  \frac
{1}{j^{\alpha}}x_{j}\right)  _{j=1}^{\infty}.
\]
Now, a characterization of Hilbert-Schmidt operators, due to Pe\l czy\'{n}ski
(see \cite{Pelc}) asserts that it suffices to show that $\Psi_{1}%
(\overset{\vee}{P})$ is a Hilbert-Schmidt operator. But is is easy to check,
since $\sum\limits_{k=1}^{\infty}\left\Vert \Psi_{1}^{(2)}(\overset{\vee}%
{P})(e_{k})\right\Vert _{l_{2}}^{2}=\sum\limits_{k=1}^{\infty}\left[  \frac
{1}{k^{\alpha}}\right]  ^{2}<\infty.$


The author wishes to acknowledge Professors M. C. Matos, G. Botelho and E.
\c{C}aliskan for important advice.

\end{document}